\documentclass[11pt, a4paper, english]{article} 
\usepackage[margin=2.5cm]{geometry}
\usepackage{babel}
\usepackage{bm}	    
\usepackage{graphicx, psfrag}
\usepackage{amsfonts, amssymb}
\usepackage{dsfont}	                        % doublestroke backage \dsfont{R} z.B. anstelle \mathbb{R}
\usepackage{yfonts} 				% Schriftgroessen fuer altdeutsche Zeichen
\usepackage{amsmath, amscd}			% AMSmath und kommutative Diagramme
\usepackage{nicefrac}	                        %boldmath griechische Symbole fett setzen

\newcommand\TRaum[2]{\ensuremath{T_{#2} #1}}							% Tangentialtraum am Punkt
\newcommand\vecf[1]{{#1}} % \bm		                        				% Vektorfeldfeld
\renewcommand\vec[1]{\bm{#1}}            					  		% Vektorfont
\newcommand\myvec{\vec}
\newcommand\LG[1]{\ensuremath{\mathit{#1}}}							% Lie-Gruppe (kursiv schreiben)
\newcommand\Lg[1]{L_{#1}}						 	                % Linkstransport
\newcommand\Rg[1]{R_{#1}}						 	                % Rechtstransport
\newcommand\liekl[2]{\left[\,#1,#2\,\right]}							% Lie-Klammer
\newcommand{\dash}[2][keinIndex]{#2'}								% Strich und 2 Strich (Ableitung nach Bogenlänge)
\newcommand{\ddash}[2][keinIndex]{#2''}

\newcommand{\pushforward}[2][keinIndex]{							% Pushforward
	\ifthenelse{\equal{keinIndex}{#1}}{\ensuremath{#2_*}}{\ensuremath{#2_{#1*}}}
	}	
	
\newcommand{\sprod}[2]{\ensuremath\left<\,#1,\, #2\,\right>}					% Skalarprodukt
\newcommand{\cprod}[2]{\ensuremath #1 \times #2}					% Skalarprodukt
\newcommand\pullback[2][keinIndex]{
	\ifthenelse{\equal{keinIndex}{#1}}{\ensuremath{#2^*}}{\ensuremath{#2^{*,#1}}}
	}

\DeclareMathOperator{\Span}{span}								% span

\def\RR{\mathbb R}										% R
\newcommand{\basis}[1]{\mathfrak{#1}}

\newcommand{\keywords}[1]{\par\noindent {\small{\em Keywords\/}: #1}}
\newtheorem{remark}{Remark}

\begin{document}\author{
					Carsten Collon{\thanks{Younicos AG, 12489 Berlin, Germany}}\ ~and Joachim Rudolph\thanks{Department of Mechatronics, Chair of Systems Theory and Control Engineering, Saarland University,\newline 66123 Saarbr\"{u}cken, Germany}}

\title{Invariant feedback control for the kinematic car on the sphere}
\maketitle

\begin{abstract}
The design of an invariant tracking control law for the kinematic car driving on a sphere is discussed. Using a Lie group framework a left-invariant description on $\LG{SO}(3)$ is derived. Basic geometric considerations allow a direct comparison of the model with the usual planar case.  Exploiting the Lie group structure an invariant tracking error is defined and a feedback is designed. Finally, one possible design of an invariant asymptotic observer is sketched.

\end{abstract}

\section{Introduction}
The model of the kinematic car, also known as unicycle, is one of the most prominent examples in nonlinear control. 
 In the present note it serves as an example for a control system for which the associated tracking control problem naturally enjoys relevant symmetry properties, i.e.\ invariance with respect to translation and rotation of the car. In general, symmetry properties are not invariant with respect to feedback. 
 This motivates the design of compatible feedback laws, denoted as invariant feedbacks, based on invariant tracking errors as proposed in \cite{RR99, Rud03b}. For a constructive approach to the computation of invariant tracking errors for systems with Lie symmetries see \cite{MRR04}. Closely related to the invariant feedback design is the design of symmetry-preserving observers, see for instance \cite{AR02, BMR08}, or \cite{BMR2009} in particular for systems on Lie groups. 

For a planar tracking problem, as for the planar motion of the car, a well-known invariant tracking error is obtained by expressing the usual tracking error with respect to a moving frame attached to the reference or vehicle trajectory, see for instance \cite{Woe1998, GR1998, RR99} for a mobile platform and the kinematic car, or \cite{Rud03b, RF03} for planar tracking for a rigid body. Other invariant tracking error candidates can be derived from geometric considerations, see for instance \cite{CRW2010} for a projection approach.
The present note considers the natural extension of the planar tracking problem to the spherical case, i.e.\ a vehicle driving on the surface of a sphere. While related problems such as the attitude control of a satellite with two controls modulo the orientation along one axis \cite{BM1999}, the control of a pointing device on $S^2$ and tracking on the sphere (without orientation) \cite{BMS1995} or the attitude control problem of the rigid-body in the fully-actuated case have been extensively considered, it seems that the spherical tracking control case for the kinematic car on a sphere has yet been overlooked, even though a possible technological application is the realization of a non-holonomic spherical joint with two controls allowing for rotation and orientation similar to the joints proposed in \cite{SNC1994}. Here, the spherical case serves as an instructive example in the discussion of the application of invariant feedback designs and its connection to structural properties of a given control problem.

The article is organized as follows.  Section~\ref{sec:lislg} recalls the basic formulation of left-invariant systems on Lie groups. In Section~\ref{sec:planar} the planar case is reviewed for later comparison with the spherical case, for which the model is derived in Section~\ref{sec:sphere}. Finally, an invariant feedback design for the spherical case and a sketch of an invariant observer is discussed in Section~\ref{sec:controlso3}.

\section{Left-invariant systems on Lie groups}
\label{sec:lislg}

In this section the class of left-invariant control systems on Lie groups is recalled. Since only the invariance property will be exploited in the following, the review is restricted to basic facts. Control systems defined on Lie groups have been initially considered by Brockett, Jurdjevic and Sussmann in the early 1970s  \cite{Bro1972, Bro1973, JS1972}, and the reader is referred to the mentioned references for a detailed discussion. For motion control and optimal control for this class of control systems see \cite{Leo1994, Jur1993, Jur1997}, a survey is given in \cite{Sac2009}. For details on Lie group theory and its application to differential equations see \cite{War83, Olv93, Boo03}.

\subsection{Lie group and Lie algebra }

An $r$-parameter Lie group is a group $G$ which is also an $r$-dimensional smooth manifold in such a way that the group multiplication the and inversion  are smooth maps. The $r$ group parameters play the role of local coordinates for $G$. 
Each element $g\in G$ defines a diffeomorphism $\Lg{g}: G \rightarrow G$, $\Lg{g}(h) = gh$ denoted as left-translation. Let $\vecf{X}$ be a vector field on $G$, i.e.\ $\vecf{X}(g)\in \TRaum{G}{g}$, and let $\pushforward{(\Lg{g})}: \TRaum{G}{}\rightarrow \TRaum{G}{}$ denote the pushforward induced by $\Lg{g}$. A vector field $\vecf{X}$ is left-invariant, if it is $\Lg{g}$-related to itself, i.e.\ if $\pushforward{(\Lg{g})}\vecf{X}(h) = \vecf{X}(gh)$ holds for all $g,h\in G$. The set of all left-invariant vector fields on $G$ forms a $\RR$-vector space $\basis{g}$ which is called the Lie algebra of $G$. Any left-invariant vector field is uniquely defined by its value at the identity $e$, thus allowing the identification $\basis{g}\simeq \TRaum{G}{e}$. Since the dimension of $\basis{g}$ equals the group dimension $r$ and all elements of the Lie algebra are left-invariant vector fields, one can choose a set of $r$ vector fields $\vecf{X}_1, \hdots,\vecf{X}_r\in \basis{g}$ with $\Span\{\vecf{X}_1(g), \hdots,\vecf{X}_r(g)\} = \TRaum{G}{g}$ for all $g\in G$.

\subsection{Left-invariant systems on Lie groups}

Using the fact that a basis for $\basis{g}$ yields also a basis for the tangent space $\TRaum{G}{g}$ smoothly depending on the base point $g$, a left-invariant system with state $g\in G$ is introduced as
\begin{align} \label{eq:lis1}
  \dot{g} &= \vecf{X}_0(g) + \sum_{k=1}^m u_k \vecf{X}_k(g),
\end{align}
with smooth left-invariant vector fields  $\vecf{X}_k$ on $G$ and smooth inputs $u_k$, $k=0, \hdots, m$. Due to the left-invariance the value of each vector field at $g$ is given by the pushforward of its value in $e$, i.e.\ $\vecf{X}_k(g)=\pushforward{(\Lg{g})} \vecf{X}_k(e)$. Consequently, the left-invariant system~\eqref{eq:lis1} can be rewritten as
\begin{align}\label{eq:lis2}
  \dot{g} &= \pushforward{(\Lg{g})}\!\left( \vecf{X}_0 + \sum_{k=1}^m u_k \vecf{X}_k\right),
\end{align}
with $\vecf{X}_k(e)=:\vecf{X}_k\in \basis{g}$, $k=0,\hdots, m$, and using the fact that the pushforward is a linear map. As usual, the vector field $\vecf{X}_0$ is called the drift vector field, and in the case $\vecf{X}_0\equiv 0$ the system is drift-free (or homogeneous \cite{JS1972}). 

In the following, the considered groups $\LG{SE}(2)$  and $\LG{SO}(3)$ have matrix representations. Consequently, tangent vectors are matrices of the same dimension, and in coordinates the pushforward of the left-translation $\Lg{g}$ is given by $g$ itself.

%----
\section{The kinematic car in the plane}
\label{sec:planar}

\begin{figure}
	\centering
  \psfrag{theta}{$\theta$}
  \psfrag{varphi}{$\varphi$}
  \psfrag{l}{$l$}
  \psfrag{v}{$v\myvec{\tau}$}
  \psfrag{V}{${\myvec{v}_F}$}
  \psfrag{y}{$\myvec{y}$}
  \psfrag{P}{$P$}
  \psfrag{pole}{\hspace{-3ex} velocity pole}
  \psfrag{osculating}{osculating circle}
  \psfrag{y1}{$y_1$}
  \psfrag{y2}{$y_2$}
	\includegraphics[width=.6\linewidth]{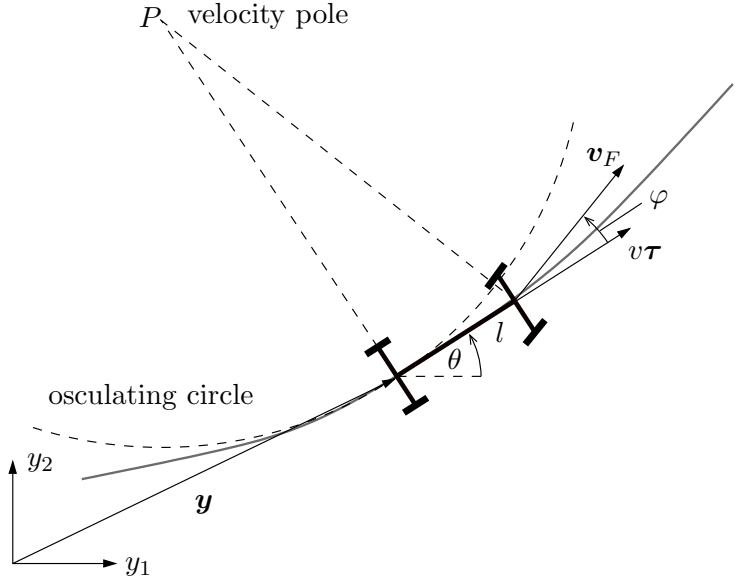}
  \caption{Kinematic car in the plane: velocity pole $P$ and osculating circle}
  \label{fig:Fzg_Modellbildung}
\end{figure}

For the sake of completeness and for later comparison the model of the kinematic car in the plane is recalled. Consider the planar motion of a car of length $l$ described by the position $\vec{y}$ of the rear axle midpoint, the steering angle $\varphi$, its orientation with respect to some inertial frame given by the angle $\theta$, and the driving speed $v$ (Figure~\ref{fig:Fzg_Modellbildung}). Assuming that the wheels roll without slipping (i.e.\ the car does not drift) one obtains the well-known model 
\begin{align}\label{eq:kinematisches_Fzg}
        \dot{\vec{y}}&= v \begin{pmatrix}\cos\theta \\ \sin\theta \end{pmatrix} = v \vec{\tau}, \quad 
	\dot{\theta} = \frac{v}{l} \tan\varphi,	
\end{align}
where $\vec{\tau}$ denotes the tangent vector to the trajectory of the rear axle midpoint. From a geometric viewpoint it seems natural to consider an arc length parametrization instead of a time parametrization. Let $s(t)=\int_0^t v(\eta) d\eta$ denote the arc length of the curve $t\mapsto \vec{y}(t)$. By using the differential relation $ds = v dt$ the system with respect to $s$ reads
\begin{align} \label{eq:fzg_ins}
  \vec{y}' &=\vec{\tau}, \quad \theta' = \frac{\tan \varphi}{l} = \kappa,
\end{align}
where $'$ denotes the derivative with respect to the arc length $s$ and $\kappa$ is the curvature of the curve. 

% Abbildung zur "Modellbildung"

% ---
\subsection{Lie group formulation}

It is well known that the model equations~\eqref{eq:kinematisches_Fzg} are form-invariant with respect to actions of elements of $\LG{SE}(2)$, i.e.\ translation and rotation, which are symmetries of the model. The symmetry property is a natural consequence of the irrelevance of the absolute position and orientation of the vehicle for its behaviour.
 In fact, the motion can be easily interpreted as a left-invariant system on $\LG{SE}(2)$. 
Hence, the motion of the car in the plane can be described by a smooth curve $\RR\supset I\owns t\mapsto g(t)\in \LG{SE}(2)$, or, equivalently, by a differential equation on the Lie group itself. Representing an element $g\in \LG{SE}(2)$ as a matrix 
\begin{align*}
   g &= \begin{pmatrix} R_{\theta} & \vec{y} \\ \vec{0} & 1 \end{pmatrix}, \quad \text{with } R_{\theta} =\begin{pmatrix} \cos \theta & -\sin\theta \\ \sin \theta & \cos\theta \end{pmatrix}\in \LG{SO}(2), 
\end{align*} 
and choosing three vector fields $\vecf{X}_i$, $i=1,2,3$ as basis for the Lie algebra $\basis{se}(2)$ 
\begin{align*}
   \basis{se}(2) = \Span \left\{\vecf{X}_1=\begin{pmatrix} 0 & -1 & 0 \\ 1 & 0 &0 \\ 0 & 0 &0 \end{pmatrix}, \, \vecf{X}_{i+1} = \begin{pmatrix} \vec{0}_{2\times 2} & \vec{e}_i \\ \vec{0}_{1\times 2} &0 \end{pmatrix}, \, i=1,2\right\},
\end{align*}
where $\vec{e}_i$ denotes the $i$-th unit vector, and following the spirit of system~\eqref{eq:lis2} the velocity of the car for the chosen nominal configuration $g=e$ is given by 
\begin{align*}
  \left.\dot{g}\right|_e &=  \begin{pmatrix} 0 & -1 & 0 \\ 1 & 0 &0 \\ 0 & 0 &0 \end{pmatrix} \frac{v}{l}\tan \varphi + \begin{pmatrix} 0 & 0 & 1 \\ 0  & 0 &0 \\ 0 & 0 &0 \end{pmatrix} v.
\end{align*}
Here, the first part describes the rotation around the origin with the angular velocity $\dot{\theta}$ and the second part models the tangential motion. Using the pushforward of the left-translation one arrives at a left-invariant model for the planar kinematic car on $\LG{SE}(2)$:
\begin{align*}
   \dot{g} &= g\left.\dot{g}\right|_e = v g   \begin{pmatrix} 0 & - \frac{\tan\varphi}{l}  & 1 \\ 
                                                                    \frac{\tan\varphi}{l} & 0 & 0 \\
                                                                    0 & 0 & 0 \end{pmatrix} = \begin{pmatrix} \frac{d}{d\theta} R_{\theta} \dot{\theta} & \vec{\tau} v \\ \vec{0} & 0 \end{pmatrix}.
\end{align*}

\subsection{Left-invariant tracking error}

Due to the Lie group structure a left-invariant tracking error is obtained by defining the usual tracking error in terms of the group multiplication. Let $t\mapsto g_d(t)$ be a smooth reference trajectory and define the tracking error by the element $g_e\in \LG{SE}(2)$ which is mapped on $g(t)$ by the desired position  $g_d(t)$, $g_e = g_d^{-1} g$. Note, that this error is also obtained as result of the normalization approach to the construction of invariant tracking errors proposed in \cite{MRR04} based on moving frames. Left-invariance with respect to any $h\in G$ is easily observed by $\tilde{g}_e = (h g_d)^{-1} (hg) =  g_e$. The construction of the left-invariant error on $\LG{SE}(2)$ yields 
\begin{align*}
   g_e &= \begin{pmatrix} R_{\theta-\theta_d} & R_{\theta_d}^T (\vec{y}-\vec{y}_d) \\ 
                          \vec{0}        & 1 \end{pmatrix} 
\end{align*}
with $g\simeq(\vec{y}, \theta)$, $g_d\simeq(\vec{y}_d, \theta_d)$. From the first entry of the second column one recovers the well-known invariant tracking error given by the usual tracking error $\vec{y}-\vec{y}_d$ parametrized with respect to the moving frame associated with the reference trajectory.

% ---
\section{The kinematic car driving on a sphere}
\label{sec:sphere}

\begin{figure}
	\centering
  \scriptsize
  \psfrag{VAGK}{\hspace{-6ex}\scriptsize \begin{parbox}[t]{12ex}{front axle \mbox{great circle}}\end{parbox}}
  \psfrag{HAGK}{\hspace{-7ex}\scriptsize\begin{parbox}[t]{12ex}{\vspace{-1ex}\flushright rear axle \mbox{great circle}}\end{parbox}}
  \psfrag{HAGK0}{\hspace{-3ex}\scriptsize\begin{parbox}[t]{14ex}{\mbox{great circle} $\theta=0$}\end{parbox}}
  \psfrag{eta}{$\theta$}
  \psfrag{gamma}{$\varphi$}
  \psfrag{pH}{$\vec{y}$}
  \psfrag{pV}{}
  \psfrag{X}{$y_1$}
  \psfrag{Y}{$y_2$}
  \psfrag{Z}{$y_3$}
  \psfrag{lRe}{}
  \psfrag{rRe}{}
  \psfrag{GKe}{}
  \psfrag{s0}{}
	\includegraphics[width=.7\linewidth]{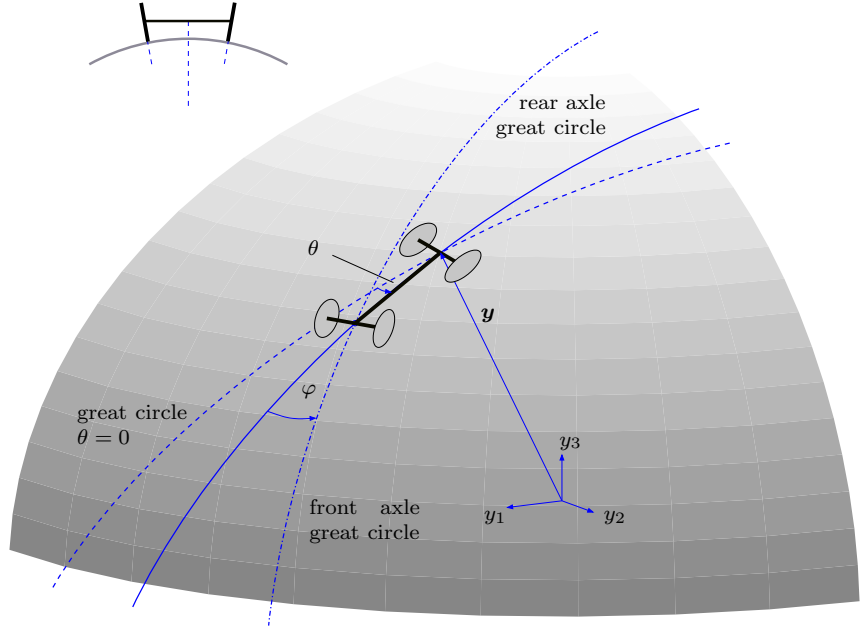}    
  \caption{Car rolling on a sphere: front view of an axle (top left), great circles corresponding to tangent and horizontal directions in the planar case }
  \label{fig:Fzg_GKs}
\end{figure}

Consider the motion of a car driving on an embedded sphere 
\begin{align*}
 S^2=\left\{(y_1,y_2,y_3)\in\RR^3\,|\, y_1^2+y_2^2+y_3^2-\rho^2=0\right\}
\end{align*}
of radius $\rho>0$. Assuming that --as shown in Figure~\ref{fig:Fzg_GKs}--  the suspension is designed in such a way that the wheels have perpendicular contact with the surface and the rolling without slipping condition still holds, at each instant of time there exist great circles through the center points of each axle corresponding to the tangent direction of motion in the planar case. Three great circles are of particular interest (Figure~\ref{fig:Fzg_GKs}): The rear axle great circle (solid line) which corresponds to the tangent to the trajectory of the rear axle center $\vec{y}$, the great circle corresponding to the horizontal through $\vec{y}$ on which the rear axle center evolves for $\theta=0$ (dashed line), and the front axle circle (dash-dotted) corresponding to the tangent to the front axle midpoint's trajectory, i.e.\ the direction of $\vec{v}_F$ in the planar case. 
As in the planar case, the two rear axle great circles enclose an angle $\theta$ describing the relative orientation of the car with respect to the nominal orientation, while the rear axle great circle and the front axle great circle enclose the steering angle $\varphi$. 

\subsection{Left-invariant model on $\LG{SO}(3)$}

\begin{figure}
  \centering  
  \scriptsize
  \psfrag{R}{$\rho$}
  \psfrag{r}{$r$}
  \psfrag{Rr}{\hspace{-3ex}$\rho+r$}
  \psfrag{L}{$l$}
  \psfrag{l}{$\ell$}
  \psfrag{lambda}{$\lambda$}
  \psfrag{HA}{rear axle}
  \psfrag{VA}{\hspace{-2ex}front axle}
  \psfrag{ge0}{$g=e$}
  \psfrag{y1}{$y_1$}
  \psfrag{y2}{$y_2$}
  \psfrag{z}{$y_3$}
  \psfrag{thetadot}{$\dot{\theta}$}
  \psfrag{vR}{$\frac{v}{\rho}$}
  \includegraphics[width=.7\linewidth]{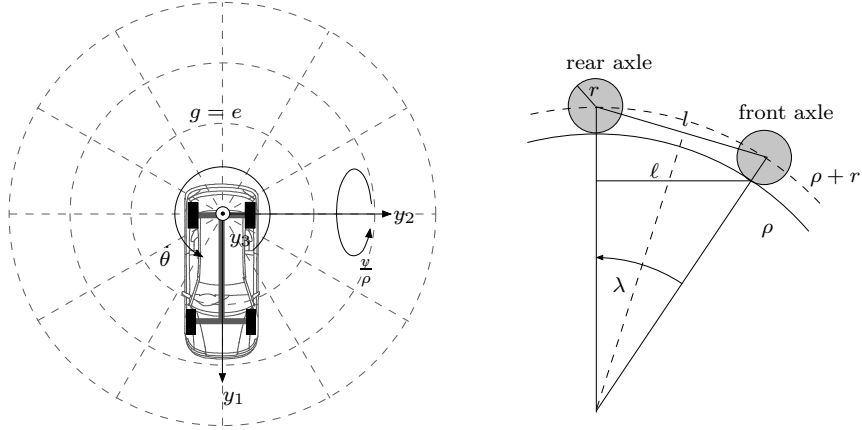}
  \caption{Left: top view of the sphere for the nominal configuration $g=e$; right: geometric construction of the length $\ell$ as funtion of the length $l$ and the sphere radius $\rho$}
   \label{fig:Rotationslaenge}
\end{figure}

The position and orientation of the car on the sphere can be uniquely described by a three-dimensional rotation matrix $g\in \LG{SO}(3)$ as follows. The matrix $g$ can be interpreted as rotated coordinate system with axes $\vec{\tau}$, $\vec{\nu}$, $\vec{\beta}$, i.e.\ $g=(\vec{\tau}, \vec{\nu}, \vec{\beta})$. Let $R(\vec{n}, \alpha)\in \LG{SO}(3)$ denote the rotation matrix associated with the rotation by $\alpha$ about the unit rotation axis $\vec{n}$. 
Associating the rotated $y_3$-axis $\vec{\beta}=g \vec{e}_3$ with the position of the rear axle center $\vec{y}$, one obtains the corresponding rotation $R_{\vec{y}}$ by setting $R_{\vec{y}} = R( \cprod{\vec{e}_3}{\frac{\vec{y}}{\rho}}, \arccos y_3 )$. (For $\vec{y}=\rho\vec{e}_3$ one sets $R_{\vec{y}}$ to the identity matrix. An alternative characterization is the rotation associated with the shortest great circle segment connecting the north pole with the point $\vec{y}$, i.e.\ the orthodrome connecting the two points.) Note that the first and the third column of $R_{\vec{y}}$ span the plane which intersects the sphere at the great circle for $\theta=0$. 
Performing a second rotation about $\frac{\vec{y}}{\rho}$ by $\theta$ yields
\begin{align*}
   g &= R\left(\frac{\vec{y}}{\rho}, \theta\right) R_{\vec{y}} =R_{\theta} R_{\vec{y}} \,.
\end{align*}
Just as for the planar case the model can be easily derived by using the Lie group approach. The Lie algebra of $\LG{SO}(3)$ is given by
\begin{align*}
  \basis{so}(3) &= \Span \left\{\vecf{X}(\vec{n})= \begin{pmatrix} 0 & -n_3 & n_2 \\ n_3 & 0 & -n_1 \\ -n_2 & n_1 &0 \end{pmatrix} \right\},
\end{align*}
where the unit vector $\vec{n}=\left(n_1,n_2,n_3\right)^T\in \RR^3$ defines the axis of rotation for the infinitesimal generator $\vecf{X}(\vec{n})$. 
Let $g=e$ denote the nominal configuration of the car, with $e=I_{3\times 3}$ being the identity matrix. In this configuration the rear axle midpoint is located at the north pole of the sphere and the car is understood to be oriented as shown on the left of Figure~\ref{fig:Rotationslaenge}. Hence, in $e$ the tangent vector of the car is given by
\begin{align} \label{eq:liso31}
  \left. \dot{g} \right|_e &= \vecf{X}(\vec{e}_2) \frac{v}{\rho} + \vecf{X}(\vec{e}_3) \dot{\theta},
\end{align}
i.e.\ the translation in $y_1$-direction at the speed $v$ corresponds to a rotation about the $y_2$-axis at the angular velocity $\frac{v}{\rho}$, and the rotation about the rear axle center corresponds to a rotation about the $y_3$-axis with the angular rate $\dot{\theta}$. However, since the surface is not planar but curved with constant curvature $\frac{1}{\rho}$, the angular velocity $\dot{\theta}$ is now obtained from $\dot{\theta} = \frac{v}{\ell} \tan \varphi$ with $\ell = \rho \sin \lambda$ and $\lambda$ denoting the central angle as depicted on the right  of Figure~\ref{fig:Rotationslaenge}. Basic geometric considerations yield the relation $\frac{\lambda}{2} = \frac{l}{2 (\rho+r)}$ with the radius $r$ of the wheels. %\footnote{Using the equality $\rho \sin \lambda = 2\rho \cos \frac{\lambda}{2} \sin \frac{\lambda}{2}$ and the series expansion for the sine it is easy to see that for $\rho\rightarrow \infty$ the planar case is obtained: $\lim_{\rho\rightarrow \infty} \rho \sin \lambda=l$.  }. 
Based on equation~\eqref{eq:liso31} a left-invariant model for the kinematic car on the sphere reads %\footnote{Again, in \cite{Jur1993} the same model but with $v\equiv 1$ is used to discuss the solution of Euler's elastica on $S^2$ by means of optimal control.}
\begin{align} \label{eq:modelso3}
  \dot{g} &= g  \left. \dot{g} \right|_e = g \left( \vecf{X}(\vec{e}_2) \frac{v}{\rho} + \vecf{X}(\vec{e}_3) \frac{v}{\ell} \tan \varphi\right) \,.
\end{align}
By construction, this model is left-invariant with respect to any element $h\in G$ corresponding to translation and rotation in the planar case.

% ----
\subsection{Geometric interpretation and differential flatness}

In the context of the previous discussion the axis  $\vec{\tau}$, $\vec{\nu}$, $\vec{\beta}$ of $g$ have a clear geometric meaning: The tangent vector $\vec{\tau}$ and the normal vector $\vec{\nu}$ span the tangent plane at $\vec{y}$ to the sphere, whereas $\vec{\beta}$ is constructed as the binormal unit vector. Further, for any given smooth trajectory $\RR\supset I \owns t \mapsto \vec{y}(t)\in \RR^3$ with $v(t)=\sprod{\dot{\vec{y}}(t)}{\dot{\vec{y}}(t)}^{\frac{1}{2}}\neq 0$ for all $t\in I$ on the embedded sphere, i.e.\ $\sprod{\vec{y}}{\vec{y}}=\rho^2$,  one has
\begin{align} \label{eq:taunuetc}
  \vec{\tau} = \frac{\dot{\vec{y}}}{v} = \vec{y}', \quad \vec{\nu}=\cprod{\frac{\vec{y}}{\rho}}{\frac{\dot{\vec{y}}}{v}}, \quad \vec{\beta} = \frac{\vec{y}}{\rho}\,.
\end{align}
Now, let $\vec{\varepsilon}_1$, $\vec{\varepsilon}_2$, $\vec{\varepsilon}_3$ denote the usual Frenet-Serret frame of the space curve 
\begin{align*}
 \vec{\varepsilon}_1 = \vec{y}' = \vec{\tau}, \quad \vec{\varepsilon}_2 = \frac{\ddash{\vec{y}}}{\kappa}, \quad \vec{\varepsilon}_3= \cprod{\vec{y}'}{\frac{\vec{y}''}{\kappa}}
\end{align*}
with $\kappa=\left\|\ddash{\vec{y}}\right\|$. In the planar case~\eqref{eq:fzg_ins} the curvature $\kappa$ of the curve $s \mapsto \vec{y}(s)$ is controlled by the steering angle $\varphi$. In the spherical case, the curvature of the space curve consists of the geodesic curvature $\kappa_g= \frac{\tan\varphi}{\ell}$ which is controlled by the steering angle and the constant curvature $\frac{1}{\rho}$ of the sphere's surface. This can be seen as follows: The columns $\vec{\tau}$,  $\vec{\nu}$, $\vec{\beta}$ of $g$ form an orthonormal frame along the space curve $s\mapsto \vec{y}(s)\in \RR^3$. Thus, the second derivative with respect to the arc length can be written as
\begin{align*}
  \ddash{\vec{y}} &= \sprod{\ddash{\vec{y}}}{\vec{\tau}} \vec{\tau}  + \sprod{\ddash{\vec{y}}}{\vec{\nu}} \vec{\nu} + \sprod{\ddash{\vec{y}}}{\vec{\beta}} \vec{\beta}.
\end{align*} 
From $\sprod{\dash{\vec{y}}}{\dash{\vec{y}}}=1$ it follows $\sprod{\ddash{\vec{y}}}{\dash{\vec{y}}}=0$, and using $\sprod{\vec{y}}{\vec{y}}=\rho^2$ one obtains $\sprod{\ddash{\vec{y}}}{{\vec{y}}} = -\sprod{\dash{\vec{y}}}{\dash{\vec{y}}}=-1$ and therefore
\begin{align*}
  \ddash{\vec{y}} &= \sprod{\ddash{\vec{y}}}{\vec{\nu}} \vec{\nu} - \frac{\vec{\beta}}{\rho} = \kappa_g \vec{\nu}- \frac{\vec{\beta}}{\rho}.
\end{align*} 
Since $\vec{\beta}$ points outward, i.e.\ $\vec{\beta}$ is the negative of the usually defined surface normal, the projection on $\vec{\beta}$ yields the negative of the surface curvature $\frac{1}{\rho}$. Comparing this with the definition of the curvature $\kappa$ of the space curve one easily obtains $\kappa^2 = \kappa_g^2+\rho^{-2}$. Now, from the model equations~\eqref{eq:modelso3} one has
\begin{align} \label{eq:kappag}
  \ddash{\vec{y}} &= \dash{\vec{\tau}} = \frac{\dot{g}}{v} \vec{e}_1 = g \begin{pmatrix} 0 & \frac{\tan\varphi}{\ell} & -\frac{1}{\rho}  \end{pmatrix}^T = \frac{\tan\varphi}{\ell} \vec{\nu} - \frac{\vec{\beta}}{\rho},
\end{align}
from which $\kappa_g$ is easily identified. For the curve on the surface $S^2$ this leads to the Frenet-Serret equations $\nabla_{\vec{\tau}}{\vec{\tau}} = \kappa_g \vec{\nu}$ and $\nabla_{\vec{\tau}}{\vec{\nu}} = -\kappa_g \vec{\tau}$, where $\nabla_{\vec{\tau}}$ denotes the covariant derivative along $\vec{\tau}$, i.e.\ the tangential part of the directional derivative (see for instance \cite{Boo03, Kue2006}). Further, by observation of equations \eqref{eq:taunuetc}, \eqref{eq:kappag}, and $v=\pm\sqrt{\vec{y}^T\vec{y}}$ it follows that the position of the rear axle center $\vec{y}$ is --just as in the planar case-- a flat output of the model.

% --- 
\section{Invariant tracking control on the sphere}
\label{sec:controlso3}

\begin{figure}
  \centering
  \scriptsize    
  \psfrag{X}{$y_1$}
  \psfrag{Y}{$y_2$}
  \psfrag{Z}{$y_3$}    
  \psfrag{FGK}{\begin{parbox}{12ex}{error circle\\ segment}\end{parbox}}
  \psfrag{TGK}{\begin{parbox}{18ex}{desired rear axle\\ great circle}\end{parbox}}
  \psfrag{delta}{$\delta$}
  \psfrag{sigma}{$\sigma$}
  \psfrag{yd}{$\vec{y}_d$}
  \psfrag{y}{$\vec{y}$}  
	\includegraphics[width=.7\linewidth]{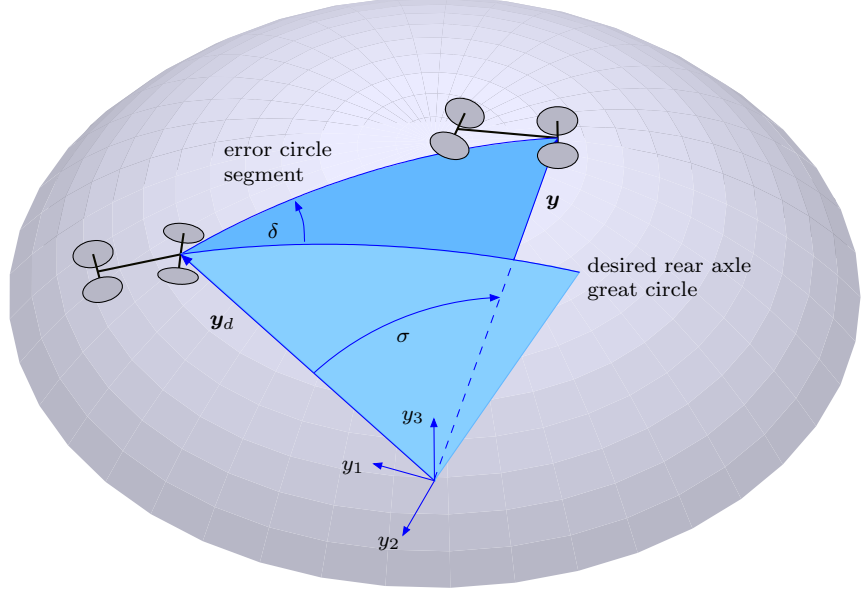}
  \caption{Parametrization of the invariant tracking error by the angles $\sigma$ and $\delta$}
  \label{fig:FolgefehlerKugel}
\end{figure}

As for the planar case, symmetry is an intrinsic property of the tracking control problem on the sphere. Hence, a feedback design is carried out using an invariant tracking approach assuming that only the position  $\vec{y}$ of the rear axle center is measured. One invariant tracking error is obtained by using a suitable parametrization of the orthodrome connecting these two points on the sphere (Figure~\ref{fig:FolgefehlerKugel}). Defining two angles $\sigma$ and $\delta$ as central angle of the orthodrome and relative orientation of the orthodrome with respect to the great circle of the rear axle on the reference trajectory, respectively.  Note that one can interpret $\sigma$ as contouring error as $\vec{y}$ moves on the desired rear axle great circle with offset $\sigma$ for $\delta=0$. (In \cite{BMS1995} the angle $\sigma$ is used as distance function on $S^2$ and the tangent direction along the error circle is denoted as ``geodesic versor''.) Similarly, $\delta$ can be interpreted as a misalignment between the current and the desired rear axle great circles. 
Using the embedding, these two angles can be determined by the relations
\begin{subequations}
\begin{align}
  \sprod{\vec{\beta}}{\vec{\beta}_d} &= \cos \sigma  \label{eq:eqsigma} \quad \text{and}  \\
  \sprod{\cprod{\vec{\beta}}{\vec{\beta}_d}}{\vec{\nu}_d} &= \cos\delta \sin \sigma  \label{eq:eqdelta}  \,.
\end{align}
\end{subequations}

\begin{remark}
 For $\vec{\beta}\neq\vec{\beta}_d$ the usual cross product ${\cprod{\vec{\beta}}{\vec{\beta}_d}}$ yields the axis of rotation of the error great circle and also defines the normal vector $\vec{\nu}$ on this great circle segment. Hence, this vector is an element of $\TRaum{S^2}{\vec{y}}$ and $\TRaum{S^2}{\vec{y}_d}$. Further, it is invariant under the pushforward along the error circle. Therefore, the projection $\sprod{\cprod{\vec{\beta}}{\vec{\beta}_d}}{\vec{\nu}_d}$ is defined on $\TRaum{S^2}{\vec{y}_d}$.
\end{remark}

In order to avoid singularities for the case $v_d=0$ it is worthwhile to use an arc length parametrization for the feedback law. Choosing the arc length $s_d(t) = \int_0^t \|\vec{y}_d(\eta)\|d\eta $ of the reference trajectory $t\mapsto \vec{y}_d$ as independent parameter one obtains the model 
\begin{align} \label{eq:arclengthso3}
  \dash{g} &= u \frac{\dot{g}}{v} = g \left(\vecf{X}(\vec{e}_2) \frac{u}{\rho} + \vecf{X}(\vec{e}_3) \frac{u}{\ell} \tan \varphi\right) = g \left.\dash{g}\right|_e, 
\end{align}
with $u=\frac{v}{v_d}$, which is considered for the following feedback design. Taking the derivative with respect to $s_d$ in equation~\eqref{eq:eqsigma} yields
\begin{align*}
  \sprod{\dash{\vec{\beta}}}{\vec{\beta}_d}\! + \! \sprod{\vec{\beta}}{\dash{\vec{\beta}}_d} &= u \sprod{\frac{\vec{\tau}}{\rho}}{\vec{\beta}_d} \! + \! \sprod{\vec{\beta}}{\frac{\vec{\tau}_d}{\rho}} = -\dash{\sigma} \sin\sigma \,
\end{align*}
and choosing an error dynamics for $\sigma$, e.g.\ $\dash{\sigma} + c_{\sigma} \sigma = 0$, $c_{\sigma}>0$, determines an invariant feedback law
\begin{align} \label{eq:feedbacku}
  u &= \rho \frac{c_{\sigma} \sin \sigma - \sprod{\vec{\beta}}{\vec{\tau}_d}}{\sprod{{\vec{\tau}}}{\vec{\beta}_d}}=\mu(\sigma, g, g_d) \,.
\end{align} 
Note that taking the limit $\vec{\beta}\rightarrow \vec{\beta}_d$, i.e.\ $\sigma\rightarrow 0$, and $\vec{\tau}\rightarrow \vec{\tau}_d$ and employing L'H\^{o}pital's rule, leads to $u \rightarrow 1$ as expected. Differentiating~\eqref{eq:eqdelta} with respect to the arc length  yields the expressions
\begin{align} \label{eq:first}
  \begin{split}
	\sprod{\cprod{\dash{\vec{\beta}}}{\vec{\beta}_d}+\cprod{{\vec{\beta}}}{\dash{\vec{\beta}}_d}}{\vec{\nu}_d} +  \sprod{\cprod{\vec{\beta}}{\vec{\beta}_d}}{\dash{\vec{\nu}}_d}   = \dash{\sigma} \cos\sigma \cos\delta - \dash{\delta} \sin \delta \sin \sigma \\
\end{split}
\end{align} 
and 
\begin{align} \label{eq:feedback2}
  \begin{split}
&\sprod{\cprod{\ddash{\vec{\beta}}}{\vec{\beta}_d} + 2\cprod{\dash{\vec{\beta}}}{\dash{\vec{\beta}}_d}+ \cprod{{\vec{\beta}}}{\ddash{\vec{\beta}}_d}}{\vec{\nu}_d}   +2\sprod{\cprod{\dash{\vec{\beta}}}{\vec{\beta}_d}+\cprod{{\vec{\beta}}}{\dash{\vec{\beta}}_d}}{\dash{\vec{\nu}}_d} +  \sprod{\cprod{\vec{\beta}}{\vec{\beta}_d}}{\ddash{\vec{\nu}}_d} \\   =&\, \ddash{\sigma} \cos\sigma \cos\delta  - (\dash{\sigma})^2 \sin \sigma \cos \delta - 2 \dash{\sigma}\dash{\delta} \sin \delta \cos \sigma 
  -\ddash{\delta} \sin \delta \sin \sigma - (\dash{\delta})^2 \cos\delta \sin \sigma \,.
  \end{split}
\end{align}
From the model equations one has
\begin{align*}
&  \dash{\vec{\beta}} = \frac{u}{\rho}\vec{\tau}, \quad \dash{\vec{\beta}}_d = \frac{\vec{\tau}_d}{\rho}, \quad \dash{\vec{\nu}}_d = -\kappa_{g,d} \vec{\tau}_d, \quad \ddash{\vec{\beta}} = \dash{u} \frac{\vec{\tau}}{\rho} + \frac{u}{\rho}\left(\kappa_g \vec{\nu} - \frac{\vec{\beta}}{\rho}\right), \\ 
&   \ddash{\vec{\beta}}_d = \frac{1}{\rho}\left(\kappa_{g,d} \vec{\nu}_d - \frac{\vec{\beta}_d}{\rho}\right), \quad \ddash{\vec{\nu}}_d = - \dash{\kappa}_{g,d} \vec{\tau}_d -\kappa_{g,d}\left(\kappa_{g,d}\vec{\nu}_d - \frac{\vec{\beta}_d}{\rho}\right)\, ,
\end{align*}
and thus, using equation~\eqref{eq:first} in order to express $\dash{\delta}$ in terms of $\sigma$, $\delta$, $\vec{\tau}$, $\vec{\beta}$, $\vec{\nu}$ and the reference trajectory, and the feedback $\mu$ from \eqref{eq:feedbacku} and its derivative $\dash{\mu}$ for $u$ and $\dash{u}$, one can fix an error dynamics for $\delta$, e.g.\ $\ddash{\delta} + c_{\delta}^1 \dash{\delta} + c^0_{\delta} \delta=0$ with positive coefficients $c_{\delta}^1$,  $c_{\delta}^0$, yielding an equation that can be locally solved for  $\kappa_g$, or, using $\kappa_g = \frac{1}{\ell} \tan\varphi$, for the steering angle $\varphi$.  Hence, an invariant control based on a quasi-static feedback can be derived. Again, taking the limit for $(\vec{\tau}, \vec{\nu}, \vec{\beta})\rightarrow (\vec{\tau}_d, \vec{\nu}_d, \vec{\beta}_d)$, i.e.\ $(\sigma, \delta) \rightarrow \vec{0}$, leads to $\kappa_g \rightarrow \kappa_{g,d}$ via  L'H\^{o}pital's rule.

% ---

\subsection{Observer design} %?

The resulting control is based on quasi-static feedback of the complete state $g\simeq(\vec{y}, \theta)$. If, as assumed in the planar case, only $\vec{y}$ is measured, an estimate for the orientation is required. This leads directly to the problem of the design of invariant asymptotic observers as introduced in \cite{AR02, AR03}. For systems defined on Lie groups a framework for invariant observer design has been presented in \cite{MBD2004, BMR2009}. However, since most observer designs on $\LG{SO}(3)$ discuss the attitude reconstruction problem for a rigid body based on gyro rate, accelerometer, and magnetometer measurements (see for instance \cite{MHP08,BMR08}), a possible local design using only the rear axle center position for the error injection is sketched at this point.

Let $\hat{g}$ denote the observer state. Consider the left-invariant observation error $g_e =g^T\hat{g}$ on $\LG{SO}(3)$.
A left-invariant observer  may be constructed as a copy of the model equations~\eqref{eq:arclengthso3} augmented by an invariant error injection, i.e.\  
\begin{align*} %\label{eq:observer1}
  \dash{\hat{g}} &= \hat{g} \left( \left. \dash{g} \right|_e + \sum_{i=1}^3 L_i(\hat{g}, \vec{y}) \vecf{X}(\vec{e}_i)   \right)
\end{align*}
with invariant observer gains $L_i$, $i=1,2,3$. From this, one obtains the differential equation for $g_e$ as
\begin{align*}% \label{eq:observer}
%\begin{split}
  \dash{g}_e &= \pushforward{(\Lg{g^T})}\dash{\hat{g}} + \pushforward{(\Rg{\hat{g}})}\dash{(g^T)}= g^T \dash{\hat{g}} - g^T \dash{g} g^T \hat{g} \\
            &= g^T \hat{g} \left(\left. \dash{g}\right|_e + \sum_{i=1}^3 L_i(\hat{g}, \vec{y})\vecf{X}(\vec{e}_i) - \hat{g}^T g  \left. \dash{g} \right|_e g^T\hat{g}  \right) \\
            &= g_e \left(\left. \dash{g}\right|_e - g_e^T \left. \dash{g} \right|_e g_e + \sum_{i=1}^3 L_i(\hat{g}, \vec{y})\vecf{X}(\vec{e}_i)  \right)\,.
%\end{split}
\end{align*}
Note that the observer is parametrized with respect to the arc length of the trajectory for $\vec{y}$. Hence, one has $u\equiv 1$ for $\left.\dash{g}\right|_e$. 
Assuming $g_e$ sufficiently close to the identity $g_e$ may be parametrized by $g_e=\exp(\vecf{X}_e)$, $\vecf{X}_e\in \basis{g}$. Using the first order approximation $g_e \simeq I_{3\times 3} + \vecf{X}_e$ one has $\vecf{X}_e \simeq  \frac{1}{2} ( g_e - g_e^T)$, i.e.\ close to $e$ the vector field $\vecf{X}_e$ is approximately the skew-symmetric component of $g_e$. Consequently, the first-order approximation of the differential equation for $\vecf{X}_e$ close to $e$ reads
\begin{align*}
  \dash{\vecf{X}}_e &=   \sum_{i=1}^3 L_i(\hat{g}, \vec{y})\vecf{X}(\vec{e}_i)  - \liekl{\left. \dash{g} \right|_e}{ \vecf{X}_e},
\end{align*}
where all higher order terms in $\vecf{X}_e$ have been neglected and $\liekl{\cdot}{\cdot}$ denotes the Lie bracket on $\LG{GL}(3, \RR)$. Since $\basis{so}(3)$ can be identified with $\RR^3$ it is possible to rewrite this matrix differential equation in vector form. For details on this identification see for instance \cite{Kod1989}. Let $\check{\vecf{X}}$ denote the unique vector $(n_1, n_2,n_3)$ with $\vecf{X}=\vecf{X}(\vec{n})$. Then, the above equation can be rewritten as
\begin{align*}
  \dash{\check{\vecf{X}}}_e &= -\left. \dash{g} \right|_e \check{\vecf{X}}_e + \sum_{i=1}^3 L_i(\hat{g}, \vec{y}) \vec{e}_i.
\end{align*}
Now, consider the left-invariant function $B=\hat{g}^T \vec{\beta}=\hat{g}^T g \vec{e}_3$. Its first-order approximation reads $B\simeq (I-(\hat{\vecf{X}}-{\vecf{X}}))\vec{e}_3 = (I- \vecf{X}_e)\vec{e}_3=\left(-\xi_2, \xi_1, 1\right)^T$, where $\check{\vecf{X}}_e=\left(\xi_1, \xi_2,\xi_3\right)^T$. Choosing the left-invariant observer gains $L_i(\hat{g}, \vec{y})=l_{i1} \sprod{B}{\hat{g}^T\hat{\vec{\nu}}}-l_{i2} \sprod{B}{\hat{g}^T\hat{\vec{\tau}}} = l_{i1} \xi_1 + l_{i2} \xi_2$ with constants $l_{ij}\in \RR$, $i=1,2,3$, $j=1,2$, one obtains the system
\begin{align*}
  \dash{\check{\vecf{X}}}_e &=  \begin{pmatrix} l_{11} & l_{12}+\kappa_g & -\frac{1}{\rho} \\
                                                     l_{21} - \kappa_g & l_{22}    & 0 \\ 
                                                     l_{31}+\frac{1}{\rho} & l_{32} & 0 \end{pmatrix}   \check{\vecf{X}}_e\, ,
\end{align*}
which is time-invariant for $l_{21}=-l_{12}=\kappa_g$. Based on this assumption, the coefficients of the corresponding characteristic polynomial
\begin{align*}
 p(\lambda) &= \lambda^3 -\left(l_{22}+ l_{11}\right) \lambda^2 +  \left(l_{31}\rho+\rho^2 l_{11}l_{22}+1\right) \frac{\lambda}{\rho^2} -\left(\frac{l_{22}}{\rho^2} + \frac{l_{31}l_{22}}{\rho}\right)
\end{align*}
can be assigned by a suitable choice of $l_{11}$, $l_{22}<0$, $l_{31}$, and hence, local asymptotic convergence of the observer error can be achieved.

% ---
\section{Conclusion}

An invariant tracking control design for the kinematic car on the sphere has been proposed. Using a Lie group framework and basic geometric considerations, the left-invariant models on $\LG{SE}(2)$ for the planar case and on  $\LG{SO}(3)$ for the spherical case can be related in a straightforward manner.


\begin{thebibliography}{10}

\bibitem{AR02}
N.~Aghannan, P.~Rouchon: On invariant asymptotic observers, in: Proc.\ 41st IEEE Conference on Decision and Control, 2002,  1479--1484.


\bibitem{AR03}
N.~Aghannan, P.~Rouchon: An intrinsic observer for a class of {Langrangian} systems, IEEE Trans.\ Autom.\ Control \textbf{48}, 2003, 936--945.

\bibitem{Boo03}
W.~M. Boothby: An introduction to differential manifolds and Riemmanian geometry, 3rd Edition, Academic Press, 2003.


\bibitem{Bro1972}
R.~W. Brockett: System theory on group manifolds and coset spaces, SIAM J.\ Control \textbf{2}, 1972,  265--284.

\bibitem{Bro1973}
R.~W. Brockett: {Lie} theory and control systems on spheres, SIAM J.\ Appl.\ Math.\ \textbf{25}, 1973,  213--225.

\bibitem{BMR08}
S.~Bonnabel, P.~Martin, P.~Rouchon: Symmetry-preserving observers, IEEE Trans.\ Autom.\ Control \textbf{53}, 2008,  2514--2526.

\bibitem{BMR2009}
S.~Bonnabel, P.~Martin, P.~Rouchon: Non-linear symmetry-preserving observers on   {Lie} groups, IEEE Trans. Autom.\ Control \textbf{54}, 2009, 1709--1713.

\bibitem{BM1999}
F.~Bullo, R.~M. Murray: Tracking for fully actuated mechanical systems: a geometric framework, Automatica \textbf{35}, 1999, 17--34.

\bibitem{BMS1995}
F.~Bullo, R.~M. Murray, A.~Sarti: Control on the sphere and reduced attitude stabilization, in: Proc.\ IFAC Symposium on Nonlinear Control Systems, 1995, 495--501.

\bibitem{CRW2010}
C.~Collon, J.~Rudolph, F.~Woittennek: {Invariant feedback design for control systems with {Lie} symmetries {--} a kinematic car example}, Discrete and Continuous Dynamical Systems, Supplement 2011, 312-–321.

\bibitem{GR1998}
D.~Guillaume, P.~Rouchon: Observation and control of a simplified car, Proc.\ IFAC Motion Control, Grenoble, 1998, 63--67.

\bibitem{Jur1993}
V.~Jurdjevic: Optimal control problems on {Lie} groups: crossroads between geometry and mechanics, in: B.~Jukbczyk, W.~Respondek (Eds.): Geometry of feedback and optimal control, Marcel Dekker, New York, 1993.

\bibitem{Jur1997}
V.~Jurdjevic: Geometric control theory, Cambridge University Press, 1997.

\bibitem{JS1972}
V.~Jurdjevic, H.~Sussmann: Control systems on {Lie} groups, J.\ Diff.\ Eq.\ \textbf{12}, 1972, 313--329.

\bibitem{Kue2006}
W.~K\"{u}hnel: Differential geometry, curves -- surfaces -- manifolds, American Mathematical Society, 2006.

\bibitem{Kod1989}
D.~E. Koditschek: The application of total energy as {Lyapunov} function for mechanical control systems, Contemporary Mathematics \textbf{97}, 1989, 131--157.

\bibitem{Leo1994}
N.~E. Leonhard: Averaging and motion control of systems on {Lie} groups, Ph.D.\ thesis, University of Maryland (1994).

\bibitem{MBD2004}
D.~Maithripala, J.~M. Berg, W.~Dayawansa: An intrinsic observer for a class of simple mechanical systems on a {Lie} group, in: Proc.\ 2004 American Control Conference, 2004.

\bibitem{MHP08}
R.~Mahony, G.~Hamel, J.-M. Pflimlin: Nonlinear complementary filters on the special orthogonal group, IEEE Trans. Autom.\ Control \textbf{53}, 2008, 1203--1218.

\bibitem{MRR04}
P.~Martin, P.~Rouchon, J.~Rudolph: Invariant tracking, ESAIM: Control, Optimisation and Calculus of Variations \textbf{10}, 2004, 1--13.

\bibitem{Olv93}
P.~Olver: Applications of Lie Groups to Differential Equations, 2nd Edition, Springer-Verlag, New-York, 1993.

\bibitem{RF03}
J.~Rudolph, R.~Fr\"{o}hlich: Invariant tracking for planar rigid body dynamics, Proc.\ Appl.\ Math.\ Mech.\ (PAMM) \textbf{2}, 2003, 9--12.

\bibitem{RR99} 
P.~Rouchon, J.~Rudolph, Invariant tracking and stabilization: problem formulation and examples, in: D.~Aeyels, F.~Lamnabhi-Lagarrigue, A.~van~der~Schaft (Eds.), Stability and stabilization of nonlinear systems, Vol.\ 246 of Lecture Notes in Control and Information Sciences, Springer-Verlag, 1999, 261--273.

\bibitem{Rud03b}
J.~Rudolph: Examples for the use of invariant errors in nonlinear control, in:  48.~Internationales Wissenschaftliches Kolloquium, Technische Universit\"{a}t  Ilmenau, 2003.

\bibitem{Sac2009}
Y.~L. Sachkov: Control theory on {Lie} groups, Journal of Mathematical Sciences \textbf{256}, 2009, 381--439.

\bibitem{SNC1994}
O.~J. S{\o}rdalen, Y.~Nakamura, W.~J. Chung: Design of a nonholonomic manipulator, in: Proc.\ 1994 IEEE Int.\ Conf.\ Robotics and Automation, 1994,
   8--13.

\bibitem{War83}
F.~W. Warner: Foundations of differentiable manifolds and Lie groups, Springer-Verlag, New-York, 1983.

\bibitem{Woe1998}
C.~Woernle: Flatness-based control of a nonholonomic mobile platform, ZAMM Suppl.\ \textbf{1}, 1998, 43--46.


\end{thebibliography}
\end{document}